\documentclass[11pt]{article}
\usepackage{graphicx}
\usepackage{amsmath,amsfonts,amssymb}
\usepackage{citesort}
\usepackage{url}
\usepackage{amsmath}
\usepackage{graphicx}
\usepackage{epsfig}
\usepackage{epstopdf}
\usepackage{color}
\def\tto{\;{\lower 1pt \hbox{$\rightarrow$}}\kern -10pt
\hbox{\raise 2pt \hbox{$\rightarrow$}}\;}

\def\ra{\rangle}
\def\la{\langle}

\def\B{I\!\!B}
\def\h{\hfill\Box}
\def\R{\Bbb R}
\def\N{\Bbb N}
\def\ox{\bar{x}}

\def\oy{\bar{y}}
\def\oz{\bar{z}}

\def\h{\hfill\square}

\def\O{\Omega}
\def\ph{\varphi}

\newcounter{lk}

\begin{document}

\begin{center}
\vspace*{0.3in} {\bf THE LOG-EXPONENTIAL SMOOTHING TECHNIQUE AND NESTEROV'S ACCELERATED GRADIENT METHOD FOR GENERALIZED SYLVESTER PROBLEMS}\\[2ex]
 N.T. An\footnote{Thua Thien Hue College of Education, 123 Nguyen Hue, Hue City, Vietnam (email: thaian2784@gmail.com).}, D. Giles\footnote{Fariborz Maseeh Department of Mathematics and Statistics, Portland State University, PO Box 751, Portland, OR 97207, United States (email: dangiles@pdx.edu). The research of Daniel Giles was partially supported by the USA National Science Foundation under grant DMS-1411817.}, N.M. Nam\footnote{Fariborz Maseeh Department of
Mathematics and Statistics, Portland State University, PO Box 751, Portland, OR 97207, United States (email: mau.nam.nguyen@pdx.edu). The research of Nguyen Mau Nam was partially supported by the USA National Science Foundation under grant DMS-1411817 and
the Simons Foundation under grant \#208785.}, R. B. Rector\footnote{Fariborz Maseeh Department of
Mathematics and Statistics, Portland State University, PO Box 751, Portland, OR 97207, United States (email: r.b.rector@pdx.edu)}.
\end{center}
\small{\bf Abstract:} The Sylvester smallest enclosing circle problem involves finding the smallest circle that encloses a finite number of points in the plane. We consider generalized versions of the Sylvester problem in which the points are replaced by sets. Based on the log-exponential smoothing technique and Nesterov's accelerated gradient method, we present an effective numerical algorithm for solving these problems.

\medskip
\vspace*{0,05in} \noindent {\bf Key words.}  log-exponential smoothing; minimization majorization algorithm; Nesterov's accelerated gradient  method; generalized Sylvester problem.

\newtheorem{Theorem}{Theorem}[section]
\newtheorem{Proposition}[Theorem]{Proposition}
\newtheorem{Remark}[Theorem]{Remark}
\newtheorem{Lemma}[Theorem]{Lemma}
\newtheorem{Corollary}[Theorem]{Corollary}
\newtheorem{Definition}[Theorem]{Definition}
\newtheorem{Example}[Theorem]{Example}
\renewcommand{\theequation}{\thesection.\arabic{equation}}
\normalsize

\section{Introduction}

The \emph{smallest enclosing circle problem} can be stated as follows: Given a finite set of points in the plane, find the circle of smallest radius that encloses all of the points. This problem was introduce in the 19th century by the English
mathematician James Joseph Sylvester (1814--1897) \cite{syl}.  It is both a facility location problem
and a major problem in computational geometry. Over a
century later, the smallest enclosing circle problem remains very active due to its important applications to clustering, nearest
neighbor search, data classification, facility
location, collision detection, computer graphics, and military operations. The problem has been widely treated in the literature from both theoretical and numerical standpoints; see \cite{AMS12,chm,ljc,FiG04,HeV81,NiN09,SVZ11,wel,Yil06,Yil07,ZTS05} and the references therein.

 The authors' recent research focuses on generalized Sylvester problems in which the given points are replaced by sets. Besides the intrinsic mathematical motivation, this question appears in  more complicated models of facility location in which the sizes of the locations are not negligible, as in bilevel transportation problems. The main goal of this paper is to develop an effective numerical algorithm for solving the \emph{smallest intersecting ball problem}.  This problems asks for the smallest ball that intersects a finite number of convex target sets in $\R^n$. Note that when the target sets given in the problem are singletons, the smallest intersecting ball problem reduces to the classical Sylvester problem.

 The smallest intersecting ball problem can be solved by minimizing a nonsmooth optimization problem in which the objective function is the maximum of the distances to the target sets. The nondifferentiability of this objective function makes it difficult to develop effective numerical algorithms for solving the problem. A natural approach is to approximate the nonsmooth objective function by a smooth function that is favorable for applying available smooth optimization schemes. Based on the log-exponential smoothing technique and Nesterov's accelerated gradient method, we present an effective numerical algorithm for solving this problem.

Our paper is organized as follows. Section 2 contains tools of convex optimization used throughout the paper. In Section 3, we focus the analysis of the log-exponential smoothing technique applied to the smallest intersecting ball problem. Section 4 is devoted to developing an effective algorithm based on the minimization majorization algorithm and Nesterov's accelerated gradient method to solve the problem. We also analyze the convergence of the algorithm.  Finally, we present some numerical examples in Section 5.

\section{Problem Formulation and Tools of Convex Optimization}
\label{s:Pre}

In this section, we introduce the mathematical models of the generalized Sylvester problems under consideration. We also present some important tools of convex optimization used throughout the paper.

Consider the linear space $\R^n$ equipped with the Euclidean norm $\|\cdot\|$. The distance function to a nonempty subset $Q$ of $\R^n$ is defined by
 \begin{equation}\label{edf}
 d(x; Q):=\inf\{\|x-q\|\; |\; q\in Q\},\; x\in \R^n.
 \end{equation}
Given $x\in \R^n$, the Euclidean projection from $x$ to is the set
$$\Pi(x; Q):=\{q\in Q\; |\; d(x; Q)=\|x-q\|\}.$$
If $Q$ is a nonempty closed convex set in $\R^n$, then $\Pi(x; Q)$ is a singleton for every $x\in \R^n$. Furthermore, the projection operator is non-expansive in the sense that
$$\|\Pi(x; Q)-\Pi(y; Q)\| \leq \|x-y\|\; \mbox{\rm for all }x, y\in \R^n.$$

Let $\O$ and $\O_i$ for $i=1,\ldots, m$ be nonempty closed convex subsets of $\R^n$. The mathematical modeling of the smallest intersecting ball problem with \emph{target sets} $\O_i$ for $i=1, \ldots, m$ and \emph{constraint set} $\O$ is
\begin{equation}\label{original_prob}
\mbox{minimize } \mathcal{D}(x):=\max\big\{d(x; \O_i)\;\big | \;i=1, \ldots, m\big\} \mbox{    subject to } x\in \O.
\end{equation}
The solution to this problem gives the center of the smallest Euclidean ball (with center in $\O$) that intersects all target sets $\O_i$ for $i=1, \ldots,m$.

In order to study new problems in which the intersecting Euclidean ball is replaced by balls generated by different norms, we consider a more general setting. Let $F$ be a closed bounded convex set that contains the origin as an interior point.  We hold this as our \emph{standing assumptions} for the set $F$ for the remainder of the paper. The \emph{support function} associated with $F$ is defined by
\begin{equation*}\label{support function}
\sigma_F(x):=\sup\{\la x, f\ra\; |\; f\in F\}.
\end{equation*}
Note that if $F=\{x\in \R^n\; |\; \|x\|_X\leq 1\}$, where $\|\cdot\|_X$ is a norm in $\R^n$, then $\sigma_F$ is the \emph{dual norm} of the norm $\|\cdot\|_X$.

Let $Q$ be a nonempty subset of $\R^n$. The generalized distance from a point $x\in \R^n$ to $Q$ generated by $F$ is given by
\begin{equation}\label{gd}
d_F(x; Q):=\inf\{\sigma_F(x-q)\; |\; q\in Q\}.
\end{equation}
The generalized distance function \eqref{gd} reduces to the distance function \eqref{edf} when $F$ is the closed unit ball of $\R^n$ with respect to the Euclidean norm. The readers are referred to \cite{nars} for important properties of the generalized distance function \eqref{gd}.

Using \eqref{gd}, a more general model of problem \eqref{original_prob} is given by
\begin{equation}\label{dfprob}
\mbox{minimize } \mathcal{D}_F(x):=\max\big\{d_F(x; \O_i)\;\big | \;i=1, \ldots, m\big\} \mbox{    subject to } x\in \O.
\end{equation}

The function $\mathcal{D}_F$ as well as its specification $\mathcal{D}$ are nonsmooth in general. Thus, problem \eqref{dfprob} and, in particular, problem \eqref{original_prob} must be studied from both theoretical and numerical view points using the tools of generalized differentiation from convex analysis.

Given a function $\ph\colon \R^n\to \R$, we say that $\ph$ is convex if  it satisfies
$$\ph(\lambda x+ (1-\lambda)y) \leq \lambda \ph(x) +(1-\lambda)\ph(y),$$
for all $x, y \in \R^n$ and $\lambda \in (0, 1)$. The function $\ph$ is said to be strictly convex if the above inequality becomes strict whenever $x\neq y$.

The class of convex functions plays an important role in many applications of mathematics, especially applications to optimization. It is well-known that for a convex function $f: \R^n \to \R$, the function has an absolute minimum on a convex set $\Omega$ at $\ox$ if and only if it has a local minimum on $\O$ at $\ox$. Moreover, if $f: \R^n\to \R$ is a differentiable convex function, then $\ox\in \O$ is a minimizer for $f$ if and only if
\begin{equation}\label{optimality}
\la \nabla f(\ox), x-\ox\ra \geq 0\; \mbox{\rm for all }x\in \O.
\end{equation}
The readers are referred to \cite{BNO03,BV04,HUL,bmn} for more complete theory of convex analysis and applications to optimization from both theoretical and numerical aspects.

\section{Smoothing Techniques for Generalized Sylvester Problems}
\label{s:SCNC}

In this section, we employ the approach developed in \cite{ZTS05} to approximate the nonsmooth optimization problem \eqref{dfprob} by a smooth optimization problem which is favorable for applying available smooth numerical algorithms. The difference here is that we use generalized distances to sets instead of distances to points.

Given an element $v\in  \R^n$, the cone generated by $v$ is given by $\mbox{\rm cone}\,\{v\}:=\{\lambda v\; |\; \lambda\geq 0\}.$ Let us review the following definition from \cite{nars}.  We recall that $F$ is a closed bounded convex set that contains zero in its interior, as per the standing assumptions in this paper.
\begin{Definition}
The set $F$ is normally smooth if for every $x\in \mbox{\rm bd}\,F$ there exists $a_x\in \R^n$ such that $N(x; F)=\mbox{\rm cone}\,\{a_x\}$.
\end{Definition}
In the theorem below, we establish the necessary and sufficient condition for the smallest intersecting ball problem \eqref{dfprob} to have a unique optimal solution.
\begin{Theorem}\label{uniqueness}
 Suppose that $F$ is normally smooth, all of the target sets are strictly convex, and at least one of the sets $\O,\O_1,...,\O_m$ is bounded. Then the smallest intersecting ball problem \eqref{dfprob} has a unique optimal solution if and only if $\bigcap_{i=1}^m \left(\O\cap \O_i\right)$ contains at most one point.
%\begin{equation}\label{original_prob}
%\mbox{minimize } \mathcal{H}(x):=\max\big\{d_F(x; \O_i)\;\big | \;i=1, \ldots, m\big\} \mbox{    subject to } x\in \O,
%\end{equation}
%has a unique solution.
\end{Theorem}
{\it Proof} It is clear that every point in the set $\bigcap_{i=1}^m \left(\O\cap \O_i\right)$ is a solution of~\eqref{dfprob}.  Thus, if~\eqref{dfprob} has a unique optimal solution we must have that $\bigcap_{i=1}^m \left(\O\cap \O_i\right)$ contains at most one point, so the necessary condition has been proven.

For the sufficient condition, assume that $\bigcap_{i=1}^m \left(\O\cap \O_i\right)$ contains at most one point.  The existence of an optimal solution is guaranteed by the assumption that at least one of the sets $\O,\O_1,...,\O_m$ is bounded.  What remains to be shown is the uniqueness of this solution.  We consider two cases.

In the first case, we assume that $\bigcap_{i=1}^m \left(\O\cap \O_i\right)$ contains exactly one point $\ox$. Observe that $\mathcal{D}_F(\ox)=0$ and $\mathcal{D}_F(x)\geq 0$ for all $x\in \R^n$, so $\ox$ is a solution of \eqref{dfprob}. If $\hat{x}\in \O$ is another solution then we must have $\mathcal{D}_F(\hat{x})=\mathcal{D}_F(\ox)=0$. Therefore, $d_F(\hat{x}; \O_i)=0$ for all $i\in \{1, \ldots, m\}$ and hence $\hat{x}\in \bigcap_{i=1}^m \left(\O\cap \O_i\right)=\{\ox\}$. We conclude that $\hat{x}=\ox$ and the problem has unique solution in this case.

For the second case we assume that $\bigcap_{i=1}^m \left(\O\cap \O_i\right)=\emptyset$.  We will show that the function
$$S(x)=\max\{\left(d_F(x; \O_1)\right)^2, \ldots, \left(d_F(x; \O_m)\right)^2\},$$
is strictly convex on $\O$. This will prove the uniqueness of the solution.

Take any $x, y\in \O$, $x\neq y$ and $t\in (0, 1)$. Denote $x_t:=tx+(1-t)y$. Let $i\in \{1, \ldots, m\}$ such that $\left(d_F(x_t; \O_i)\right)^2=S(x_t)$.  Let $u, v \in \O_i$ such that $\sigma_F(x-u)=d_F(x; \O_i)$ and $\sigma_F(y-v)=d_F(y; \O_i)$. Then we have
\begin{align*}
S(x_t)&=\left(d_F(x_t; \O_i)\right)^2\\
&=\left[d_F(tx+(1-t)y; \O_i)\right]^2\\
&\leq \left[  td_F(x; \O_i) + (1-t)d_F(y; \O_i)\right]^2\\
&= \left[  t\sigma_F(x-u) + (1-t)\sigma_F(y-v)\right]^2\\
&= t^2\left(\sigma_F(x-u)\right)^2 +2t(1-t)\sigma_F(x-u)\sigma_F(y-v) + (1-t)^2\left(\sigma_F(y-v)\right)^2\\
&\leq  t^2\left(\sigma_F(x-u)\right)^2 +t(1-t)\left[\left(\sigma_F(x-u)\right)^2 + \left(\sigma_F(y-v)^2\right)\right] + (1-t)^2\left(\sigma_F(y-v)\right)^2\\
&=t\left(\sigma_F(x-u)\right)^2 + (1-t)\left(\sigma_F(y-v)\right)^2 \\
&= t\left(d_F(x; \O_i)\right)^2 +  (1-t)\left(d_F(y; \O_i)\right)^2\\
&\leq tS(x) +  (1-t)S(y).
\end{align*}
Recall that we need to prove the inequality $S(x_t)<tS(x)+(1-t)S(y)$.  Suppose by contradiction that $S(x_t)=tS(x)+(1-t)S(y)$. Then all of the inequalities in the above estimate are turned to equalities and thus we have
\begin{equation*}
d_F(x_t; \O_i)=td_F(x; \O_i)+ (1-t)d_F(y; \O_i)
\end{equation*}
and
\begin{equation}\label{est_uniq_0}
\sigma_F(x-u)=\sigma_F(y-v).
\end{equation}
Hence,
\begin{equation}
d_F(x_t; \O_i)=\sigma_F(x-u)=\sigma_F(y-v).
\label{est_uniq}
\end{equation}
Observe that $\sigma_F(w)=0$ if and only if $w=0$, \eqref{est_uniq} implies $x=u$ if and only if $y=v$. We claim that $x\neq u$ and $y\neq v$. Indeed, if $x=u$ and $y=v$, then $x, y \in \O_i$ and hence $x_t\in \O_i$ by the convexity of $\O_i$. Thus $d_F(x_t; \O_i)=0$. This contradicts the fact that  $d_F(x_t; \O_i)=\mathcal{D}_F(x_t)>0$  which guaranteed by the assumption $\bigcap_{i=1}^m \left(\O \cap \O_i \right) = \emptyset$.

Now, we will show that $u\neq v$. Denote $c:=tu+(1-t)v\in\O_i$. The properties of the support function and \eqref{est_uniq_0} give
\begin{align*}
d_F(x_t; \O_i)&\leq \sigma_F\left(x_t - c\right) \\
&=\sigma_F\left(t(x-u)+(1-t)(y-v)\right)\\
& \leq \sigma_F\left(t(x-u)\right)+\sigma_F\left((1-t)(y-v)\right)\\
&\leq  t\sigma_F(x-u)+ (1-t) \sigma_F(y-v)\\
&=\sigma_F(x-u).
\end{align*}
By \eqref{est_uniq} we have $\sigma_F\left(t(x-u)+(1-t)(y-v)\right) =\sigma_F\left(t(x-u)\right)+\sigma_F\left((1-t)(y-v)\right).$  Since $F$ is normally smooth, it follows from \cite[Remark 3.4]{nars} that there exists $\lambda>0$ satisfying $$t(x-u)=\lambda(1-t)(y-v).$$
Now, by contradiction, suppose $u=v$. Then $x-u=\beta(y-u)$ where $\beta=t^{-1}\lambda(1-t)$. Note that $\beta\neq 1$ since $x\neq y$ and $\sigma_F(y-u)>0$ since $y-u\neq 0$. Now we have
\begin{align*}
\sigma_F(x-u)&=\sup\{\la x-u, f\ra \;|\; f\in F\}=\sup\{\la \beta(y-u), f\ra \;|\; f\in F\}\\
&=\beta\sup\{\la y-u, f\ra \;|\; f\in F\}=\beta \sigma_F(y-u) \neq \sigma_F(y-u),
\end{align*}
which contradicts \eqref{est_uniq}. Thus $u\neq v$.

Since $u, v\in \O_i$, $u\neq v$ and $\O_i$ is strictly convex, $c\in \mbox{int }\O_i$. The assumption $\bigcap_{i=1}^m \left(\O \cap \O_i \right) = \emptyset$ gives $d_F(x_t; \O_i)=\mathcal{D}_F(x_t)>0$.  Therefore, $x_t\notin \O_i$ and thus $x_t\neq c$. Let $\delta>0$ such that $\B(c; \delta)\subset \O_i$. Then $c+ \gamma (x_t-c)\in \O_i$, with $\gamma=\frac{\delta}{2\|x_t-c\|}>0$. We have
\begin{align*}
d_F(x_t; \O_i)&\leq \sigma_F\left(x_t - c - \gamma (x_t-c)\right) =(1-\gamma)\sigma_F\left(t(x-u)+(1-t)(y-v)\right)\\
&< \sigma_F\left(t(x-u)+(1-t)(y-v)\right)=\sigma_F(x-u).
\end{align*}
This contradicts \eqref{est_uniq} and completes the proof. $\h$

Recall the following definition.
\begin{Definition}
A convex set $F$ is said to be normally round if $N(x; F)\neq N(y; F)$ whenever $x,y\in \mbox{\rm bd}\,F$, $x\neq y$.
\end{Definition}

\begin{Proposition}\label{differentiability} Let $\Theta$ be a nonempty closed convex subset of $\R^n$. Suppose that $F$ is normally smooth and normally round. Then the function $g(x):=[d_F(x; \Theta)]^2$, $x\in \R^n$, is continuously differentiable.
\end{Proposition}
{\it Proof} It suffices to show that $\partial g(\ox)$ is a singleton for every $\ox\in \R^n$. By \cite{bmn}, we have
\begin{equation*}
\partial g(\ox)=2d_F(\ox; \Theta)\partial d_F(\ox; \Theta).
\end{equation*}
 It follows from \cite[Proposition 4.3 (iii)]{nars} that $g$ is continuously differentiable on $\Theta^c$, and so
\begin{equation*}
\partial g(\ox)=2d_F(\ox; \Theta)\nabla d_F(\ox; \Theta)=2d_F(\ox; \Theta)\nabla \sigma_F(\ox-w),
\end{equation*}
where $w:=\Pi_F(\ox; \Theta)$ and $\ox\notin \Theta$.

In the case where $\ox\in \Theta$, one has $d_F(\ox; \Theta)=0$, and hence
\begin{equation*}
\partial g(\ox)=2d_F(\ox; \Theta)\partial d_F(\ox; \Theta)=\{0\}.
\end{equation*}
The proof is now complete.
$\h$

If all of the target sets have a common point which belongs to the constraint set, then such a point is a solution of problem (\ref{dfprob}), so we always assume that $\bigcap_{i=1}^n \left(\O_i \cap \O\right) = \emptyset$. We also assume that at least one of the sets $\O,\O_1,...,\O_m$ is bounded which guarantees the existence of an optimal solution; see \cite{nh}.  These are our \emph{standing assumptions} for the remainder of this section.

Let us start with some useful well-known results.  We include the proofs for the convenience of the reader.
\begin{Lemma}\label{lm1} Given positive numbers $a_i$ for $i=1, \ldots,m$, $m>1$, and $0<s<t$, one has\\[1ex]
{\rm\bf (i)} $(a_1^s+a_2^s+\ldots +a_m^s)^{1/s}>(a_1^t+a_2^t+\ldots +a_m^t)^{1/t}$.\\
{\rm\bf (ii)} $(a_1^{1/s}+a_2^{1/s}+\ldots +a_m^{1/s})^{s}<(a_1^{1/t}+a_2^{1/t}+\ldots +a_m^{1/t})^{t}$.\\
{\rm\bf (iii)} $\lim_{r\to 0^+} (a_1^{1/r}+a_2^{1/r}+\ldots +a_m^{1/r})^{r}=\max\{a_1, \ldots, a_m\}$.
\end{Lemma}
{\it Proof} {\bf (i)} Since $t/s>1$, it is obvious that
\begin{equation*}
\Big (\dfrac{a_1^s}{\sum_{i=1}^m a_i^s}\Big)^{t/s}+\cdots +\Big(\dfrac{a_m^s}{\sum_{i=1}^m a_i^s}\Big)^{t/s}< \dfrac{a_1^s}{\sum_{i=1}^m a_i^s}+\cdots+\dfrac{a_m^s}{\sum_{i=1}^m a_i^s}=1,
\end{equation*}
since  $\dfrac{a_i^s}{\sum_{i=1}^m a_i^s}\in (0,1)$. It follows that
\begin{equation*}
\dfrac{a_1^t}{(\sum_{i=1}^m a_i^s)^{t/s}}+\cdots+\dfrac{a_m^t}{(\sum_{i=1}^m a_i^s)^{t/s}}<1,
\end{equation*}
and hence
\begin{equation*}
a_1^t+\cdots +a_m^t<(\sum_{i=1}^m a_i^s)^{t/s}.
\end{equation*}
This implies {\bf(i)} by rasing both sides to the power of $1/t$. \\[1ex]
{\bf (ii)} Inequality {\bf (ii)} follows directly from {\bf (i)}.\\[1ex]
{\bf (iii)} Defining $a:=\max\{a_1, \ldots, a_m\}$ yields
\begin{equation*}
a\leq (a_1^{1/r}+a_2^{1/r}+\ldots +a_m^{1/r})^{r}\leq m^r a \to a \;\mbox{\rm as } r\to 0^+,
\end{equation*}
which implies {\bf (iii)} and completes the proof.
$\h$

For $p>0$ and for $x\in \R^n$, the \emph{log-exponential smoothing function} of $\mathcal{D}_F(x)$ is defined  as \begin{equation}\label{dfsmooth_approx}
\mathcal{D}_F(x, p):=p \ln\sum_{i=1}^m \exp\left(\frac{G_{F,i}(x, p)}{p}\right),
\end{equation}
where
$$G_{F,i}(x, p):=\sqrt{d_F(x; \O_i)^2 +p^2}.$$
\begin{Theorem}\label{lm3.1} The function $\mathcal{D}_F(x, p)$ defined in {\rm(\ref{dfsmooth_approx})} has the following properties:\\[1ex]
{\textnormal{\bf(i)}} If $x\in \R^n$ and $0<p_1<p_2$, then
$$\mathcal{D}_F(x, p_1)<\mathcal{D}_F(x, p_2).$$
{\textnormal{\bf(ii)}} For any $x\in \R^n$ and $p>0$,
$$0\leq \mathcal{D}_F(x, p)-\mathcal{D}_F(x)\leq p(1+\ln m).$$
{\textnormal{\bf(iii)}} For any $p>0$, the function $\mathcal{D}_F(\cdot, p)$ is convex. If we suppose further that $F$ is normally smooth and the sets $\O_i$ for $i=1, \ldots, m$ are strictly convex and not collinear (i.e., it is impossible to draw a straight line that intersects all the sets $\O_i$), then $\mathcal{D}_F(\cdot, p)$ is strictly convex. \\[1ex]
{\textnormal{\bf(iv)}} For any $p>0$, if $F$ is normally smooth and normally round, then $\mathcal{D}_F(\cdot, p)$ is continuously differentiable.\\[1ex]
{\textnormal{\bf(v)}} If at least one of the target sets $\O_i$ for $i=1,\ldots,m$ is bounded, then $\mathcal{D}_F(\cdot, p)$ is coercive in the sense that
\begin{equation*}
\lim_{\|x\|\to \infty}\mathcal{D}_F(x, p)=\infty.
\end{equation*}
\end{Theorem}
{\it Proof} {\bf (i)} Define
\begin{align*}
&a_i(x,p):=\exp(G_{F,i}(x,p)),\\
&a_\infty(x,p):=\max_{i=1,\ldots,m}a_i(x,p),\;\text{and}\\
&G_{F,\infty}(x,p):=\max_{i=1,\ldots,m}G_{F,i}(x,p).
\end{align*}
Then $a_i(x,p)$ is strictly increasing on $(0,\infty)$ as a function of $p$ and
\begin{align*}
&\dfrac{a_i(x,p)}{a_\infty(x,p)}=\exp(G_{F,i}(x,p)-G_{F,\infty}(x,p))\leq 1,\\
&G_{F,\infty}(x,p)\leq \mathcal{D}_F(x)+p.
\end{align*}
For $0<p_1<p_2$, it follows from Lemma \ref{lm1}(ii) that
\begin{align*}
\mathcal{D}_F(x,p_1)&=\ln \big[\sum_{i=1}^m(a_i(x,p_1))^{1/p_1}\big]^{p_1}< \ln \big[\sum_{i=1}^m(a_i(x,p_1))^{1/p_2}\big]^{p_2}\\
&< \ln\big[\sum_{i=1}^m(a_i(x,p_2))^{1/p_2}\big]^{p_2}=\mathcal{D}_F(x,p_2),
\end{align*}
which justifies {\bf (i)}.\\[1ex]
{\bf (ii)} %Given any $r$ and $p$ such that $0<r<p$, one has $\mathcal{D}_F(x,r)<\mathcal{D}_F(x,p)$. Letting $r\to 0^+$ and following the proof of Lemma \ref{lm1}(iii) yields $\mathcal{D}_F(x)=\lim_{r\to 0^+}\mathcal{D}_F(x,r)\leq \mathcal{D}_F(x,p)$.
It follows from \eqref{dfsmooth_approx} that for any $i\in \{1, \ldots, m\}$, we have
$$\mathcal{D}_F(x,p)\geq p\ln \exp\left(\dfrac{G_{F, i}(x, p)}{p}\right)=G_{F, i}(x, p)\geq d_F(x; \O_i).$$
This implies $\mathcal{D}_F(x,p)\geq \mathcal{D}_F(x)$ for all $x\in \R^n$ and $p>0$. Moreover,
\begin{align*}
\mathcal{D}_F(x,p)=&\ln a_\infty(x,p)\Big[\sum_{i=1}^m \big(\dfrac{a_i(x,p)}{a_\infty(x,p)}\big)^{1/p}\Big]^{p}\\
&=\ln a_\infty(x,p)+p\ln\sum_{i=1}^m \big(\dfrac{a_i(x,p)}{a_\infty(x,p)}\big)^{1/p}\\
&\leq G_{F,\infty}(x,p)+p\ln m\leq \mathcal{D}_F(x)+p+p\ln m.
\end{align*}
Thus, {\bf (ii)} has been proved.\\[1ex]
{\bf (iii)} Given $p>0$, the function $f_p(t):=\frac{\sqrt{t^2+p^2}}{p}$ is increasing and convex on the interval $[0, \infty)$, and  $d(\cdot; \O_i)$ is convex, so the function $k_i(x, p):=\frac{G_{F, i}(x, p)}{p}$ is also convex with respect to $x$. For any $x, y \in \R^n$ and $\lambda \in (0, 1)$, by the convexity of the function $$u=(u_1, \ldots, u_m) \to \ln\sum_{i=1}^m  \exp(u_i),$$ one has
\begin{align}
\mathcal{D}_F(\lambda x + (1-\lambda)y, p) &= p\ln\sum_{i=1}^m \exp\bigg [k_i(\lambda x + (1-\lambda) y, p)\bigg ] \notag \\
&\leq p\ln\sum_{i=1}^m \exp\bigg [\lambda k_i(x, p) + (1-\lambda)k_i(y, p)\bigg]  \label{est8} \\
&\leq \lambda p\ln\sum_{i=1}^m \exp\bigg (k_i(x, p)\bigg)  + (1-\lambda)p\ln\sum_{i=1}^m \exp\bigg (k_i(y, p)\bigg) \notag\\
&=\lambda \mathcal{D}_F(x, p) + (1-\lambda) \mathcal{D}_F(y, p).\notag
\end{align}
Thus, $\mathcal{D}_F(\cdot, p)$ is convex. Suppose that $F$ is normally smooth and the sets $\O_i$ for $i=1, \ldots, m$ are strictly convex and not collinear, but $\mathcal{D}_F(\cdot, p)$ is not strictly convex. Then there exist $x, y\in \R^n$ with $x\neq y$ and $0<\lambda <1$ such that
$$\mathcal{D}_F(\lambda x +(1-\lambda)y, p) = \lambda\mathcal{D}_F(x, p)+ (1-\lambda)\mathcal{D}_F(y, p).$$
Thus, all the inequalities \eqref{est8} become equalities. Since the functions $\ln, \exp$ are  strictly increasing on $(0, \infty)$, this implies
\begin{equation}
k_i(\lambda x +(1-\lambda)y, p)=\lambda k_i(x, p) +(1-\lambda)k_i(y, p), \qquad \mbox{for all } i=1, \ldots, m.
\label{est9}
\end{equation}
Observe that $k_i(\cdot, p)=f_p\left(d_F(\cdot; \O_i)\right)$ and the function $f_p(\cdot)$ is strictly increasing $[0, \infty)$, it follows from \eqref{est9} that
$$d_F(\lambda x +(1-\lambda y), \O_i)=\lambda d_F(x, \O_i) +(1-\lambda)d_F(y, \O_i) \qquad \mbox{for all } i=1, \ldots, m.$$
The result now follows directly from the proof of \cite[Proposition 4.5]{nars}. \\[1ex]
{\bf (iv)} Let $\ph_i(x):=[d_F(x; \O_i)]^2$. Then $\ph_i$ is continuously differentiable by Proposition \ref{differentiability}.   By the chain rule, for any $p>0$, the function $\mathcal{D}_F(x, p)$ is continuously differentiable as a function of $x$.\\[1ex]
 {\bf (v)} Without loss of generality, we assume that $\O_1$ is bounded. It then follows from {\bf(ii)} that
$$\lim \limits_{\|x\|\to \infty}\mathcal{D}_F(x, p)\geq \lim \limits_{\|x\|\to \infty}\mathcal{D}_F(x)\geq \lim \limits_{\|x\|\to \infty} d_F(x; \O_1) =\infty.$$
Therefore, $\mathcal{D}_F(\cdot, p)$ is coercive, which justifies {\bf(iv)}. The proof is now complete.
$\h$

In the next corollary, we obtain an explicit formula of the gradient of the log-exponential approximation of $\mathcal{D}$ in the case where $F$ is the closed unit ball of $\R^n$. For $p>0$ and for $x\in \R^n$, define
\begin{equation}\label{smooth_approx}
\mathcal{D}(x, p):=p \ln\sum_{i=1}^m \exp\left(\frac{G_{i}(x, p)}{p}\right),
\end{equation}
where
$$G_{i}(x, p):=\sqrt{d(x; \O_i)^2 +p^2}.$$

\begin{Corollary} For any $p>0$, $\mathcal{D}(\cdot, p)$ is continuously differentiable with the gradient in $x$ computed by
\begin{equation*} \label{nabla}
\nabla_x \mathcal{D}(x, p) = \sum_{i=1}^m \frac{\Lambda_i(x, p)}{G_i(x, p)}\left(x-\widetilde{x}_i\right),
\end{equation*}
where $\widetilde{x}_i:=\Pi(x; \O_i)$, and
$$\Lambda_i(x, p): = \frac{\exp\left(G_i(x, p)/p\right)}{\sum_{i=1}^m\exp\left(G_i(x, p)/p\right)}.$$
\end{Corollary}
{\it Proof}  It follows from Theorem \ref{lm3.1} that $\mathcal{D}(\cdot, p)$ is continuously differentiable. Let $\ph_i(x):=[d(x; \O_i)]^2$. Then $\nabla \ph_i(x) = 2(x-\widetilde{x}_i),$ where $\widetilde{x}_i:=\Pi(x; \O_i)$, and hence the gradient formula for $\mathcal{D}(x,p)$ follows from the chain rule.
$\h$

\begin{Remark}
{\rm
{\bf (i)} To avoid working with large numbers when implementing algorithms for \eqref{original_prob}, we often use the identity
$$\Lambda_i(x, p): = \frac{\exp\left(G_i(x, p)/p\right)}{\sum_{i=1}^m\exp\left[(G_i(x, p)/p\right)}=\frac{\exp\left[(G_i(x, p)-G_\infty(x,p)]/p\right)}{\sum_{i=1}^m\exp\left[(G_i(x, p)-G_\infty(x,p)]/p\right)},$$
where $G_{\infty}(x,p):=\max_{i=1,\ldots,m}G_{i}(x,p)$.\\[1ex]
{\bf (ii)} In general, $\mathcal{D}(\cdot, p)$ is not strictly convex. For example, in $\R^2$, consider the sets $\O_1=\{-1\}\times [-1, 1]$ and $\O_2=\{1\}\times [-1, 1]$. Then $\mathcal{D}(\cdot, p)$ takes constant value on $\{0\}\times [-1, 1]$.
}\end{Remark}
An important relation between problem {\rm(\ref{dfprob})} and problem of minimizing the function (\ref{smooth_approx}) on $\O$ is given in the proposition below. Note that the assumption of the proposition involves the uniqueness of an optimal solution to problem {\rm(\ref{dfprob})} which is guaranteed under our standing assumptions by Theorem \ref{uniqueness}.
\begin{Proposition}
Let $\{p_k\}$ be a sequence of positive real numbers converging to 0. For each $k$, let $y_k \in \mbox{\rm arg}\min_{x\in \O} \mathcal{D}_F(x, p_k)$. Then $\{y_k\}$ is a bounded sequence and every subsequential limit of $\{y_k\}$ is an optimal solution of problem {\rm(\ref{dfprob})}. Suppose further that problem {\rm(\ref{dfprob})} has a unique optimal solution. Then $\{y_k\}$ converges to that optimal solution.
\end{Proposition}
{\it Proof}  First, observe that $\{y_k\}$ is well defined because of the assumption that at least one of the sets $\O,\O_1,...,\O_m$ is bounded and the coercivity of $\mathcal{D}_F(\cdot, p_k)$.  By Theorem \ref{lm3.1} (ii), for all $x\in \O$, we have
\begin{equation*}\label{ee}\mathcal{D}_F(x, p_k) \leq \mathcal{D}_F(x) + p_k(1+ \ln m)\;\; \mbox{\rm and }\;\;\mathcal{D}_F(y_k) \leq \mathcal{D}_F(y_k, p_k) \leq  \mathcal{D}_F(x, p_k).
\end{equation*}
Thus, $\mathcal{D}_F(y_k)\leq \mathcal{D}_F(x) + p_k(1+ \ln m)$, which implies the bounded property of $\{y_k\}$
using the boundedness of $\O$ or the coercivity of $\mathcal{D}_F(\cdot)$ from Theorem \ref{lm3.1} (v). Suppose that the subsequence $\{y_{k_l}\}$ converges to $y_0$. Then $\mathcal{D}_F(y_0) \leq \mathcal{D}_F(x)$ for all $x\in \O$, and hence $y_0$ is an optimal solution of problem {\rm(\ref{dfprob})}. If {\rm(\ref{dfprob})} has a unique optimal solution $\oy$, then $y_0=\oy$ and hence $y_k\to\oy$.
$\h$

Recall that a function $\ph:Q \to\R$ is called {\em strongly convex} with modulus $m>0$ on a convex set $Q$ if $\ph(x) - \frac{m}{2}\|x\|^2$ is a convex function on $Q$. From the definition, it is obvious that any strongly convex function is also strictly convex. Moreover, when $\ph$ is twice differentiable, $\ph$ is strongly convex with modulus $m$ on an open convex set $Q$ if $\nabla^2 \ph(x) - mI $ is positive semidefinite for all $x\in Q$; see \cite[Theorem 4.3.1(iii)]{HUL}.

\begin{Proposition}\label{lemma2}
Suppose that all the sets $\O_i$ for $i=1,\ldots,m$ reduce to singletons. Then for any $p>0$, the function $\mathcal{D}(\cdot, p)$ is strongly convex on any bounded convex set, and $\nabla_x \mathcal{D}(\cdot,p)$ is globally Lipschitz continuous on $\R^n$ with Lipschitz constant $\frac{2}{p}$.
\end{Proposition}
{\it Proof} Suppose that $\O_i=\{c_i\}$ for $i=1,\ldots,m$. Then
\begin{equation*}
\mathcal{D}(x, p)=p \ln\sum_{i=1}^m \exp\left(\frac{g_i(x, p)}{p}\right),
\end{equation*}
and the gradient of $\mathcal{D}(\cdot, p)$ at $x$ becomes
\begin{equation*}
\nabla_x \mathcal{D}(x, p) = \sum_{i=1}^m \frac{\lambda_i(x, p)}{g_i(x, p)}\left(x-c_i\right),
\end{equation*}
where
\begin{equation*}
g_i(x, p):=\sqrt{\|x-c_i \|^2+p^2} \; \; \mbox{\rm and }\; \lambda_i(x, p): = \frac{\exp\left(g_i(x, p)/p\right)}{\sum_{i=1}^m\exp\left(g_i(x, p)/p\right)}.
\end{equation*}
Let us denote
$$Q_{ij}:=\frac{(x-c_i)(x-c_j)^T}{g_i(x, p)g_j(x, p)}.$$
Then
$$\nabla_x^2\mathcal{D}(x, p)=\sum_{i=1}^m\left[\frac{\lambda_i(x, p)}{g_i(x, p)}(I_n-Q_{ii})+\frac{\lambda_i(x, p)}{p}Q_{ii}-\sum_{j=1}^m\frac{\lambda_i(x, p)\lambda_j(x, p)}{p}Q_{ij}\right].$$
Given a positive constant $K$, for any $x\in \R^n$, $\|x\| <  K$ and $z\in \R^n$, $z\neq 0$, one has
\begin{align*}
\dfrac{1}{g_i(x,p)}(\|z\|^2-z^TQ_{ii}z) &\geq \dfrac{1}{g_i(x,p)}(\|z\|^2 - \|z\|^2\|{(x-c_i)}/{g_i(x, p)}\|^2)\\
&=\dfrac{1}{\sqrt{\|x-c_i\|^2+p^2}} \|z\|^2\left[\frac{p^2}{\|x-c_i\|^2+p^2}\right]\\
&\geq  \|z\|^2\left[\frac{p^2}{[2(\|x\|^2+\|c_i\|^2)+p^2]^{3/2}}\right] \\
&\geq \ell \|z\|^2,
\end{align*}
where $$\ell: = \frac{p^2}{[2K+2\max_{1\leq i\leq m}\|c_i\|^2 + p^2]^{3/2}}.$$
For $m$ real numbers $a_1, \ldots, a_m$, since $\lambda_i(x, p) \geq 0$ for all $i=1, \ldots, m$, and $\sum_{i=1}^m\lambda_i(x, p)=1$, by Cauchy-Schwartz inequality, we have
\begin{equation*}
\left(\sum_{i=1}^m \lambda_i(x, p)a_i\right)^2=\left(\sum_{i=1}^m \sqrt{\lambda_i(x, p)}\sqrt{\lambda_i(x, p)}\,a_i\right)^2\leq \sum_{i=1}^m \lambda_i(x, p)a_i^2.
\end{equation*}
This implies
\begin{align*}
z^T\nabla_x^2\mathcal{D}(x, p)z &= \sum_{i=1}^m\left[\frac{\lambda_i(x, p)}{g_i(x, p)}(\|z\|^2-z^TQ_{ii}z)\right]\\
&+ \sum_{i=1}^m \left[\frac{\lambda_i(x, p)}{p}z^TQ_{ii}z-\sum_{j=1}^m\frac{\lambda_i(x, p)\lambda_j(x, p)}{p} z^TQ_{ij}z\right].\\
&\geq  \ell \|z\|^2 + \frac{1}{p}\left[\sum_{i=1}^m \lambda_i(x, p)a_i^2 - \left(\sum_{i=1}^m \lambda_i(x, p)a_i\right)^2\right] \\
&\geq  \ell \|z\|^2,
\end{align*}
where $a_i:=z^T(x-c_i)/ g_i(x, p)$. This shows that $\mathcal{D}(x, p)$ is strongly convex on $B(0; K)$.

The fact that for any $p>0$, the gradient of $\mathcal{D}(x, p)$ with respect to $x$ is Lipschitz continuous with constant $L=\frac{2}{p}$ was proved in \cite[Proposition~2]{Z07}.
$\h$

\section{The Minimization Majorization Algorithm for Generalized Sylvester Problems}

In this section, we apply the minimization majorization well known in computational statistics along with the log-exponential smoothing technique developed in the previous section to develop an algorithm for solving the smallest intersecting ball problem. We also provide some examples showing that minimizing functions that involve distances to convex sets not only allows to study a generalized version of the smallest enclosing circle problem, but also opens up the possibility of applications to other problems of constrained optimization.

Let $f: \R^n \to \R$ be a convex function. Consider the optimization problem
\begin{equation}\label{general problem}
\mbox{\rm minimize }f(x)\; \mbox{\rm subject to } x\in \O.
\end{equation}
A function $g: \R^n\to \R$ is called a \emph{surrogate} of $f$ at $\bar z \in \O$ if
\begin{align*}
& f(x)\leq g(x)\; \mbox{\rm for all }x\in \O,\\
& f(\bar z) =g(\bar z).
\end{align*}
The set of all surrogates of $f$ at $\oz$ is denoted by $S(f,\oz)$.

The minimization majorization algorithm for solving \eqref{general problem} is given as follows; see \cite{Mairal}.
\newpage

{\small{\bf Algorithm 1}.
\begin{center}
\begin{tabular}{| l |}
\hline
{\small INPUT}: $x_0\in \O$, $N$\\
{\bf for} $k=1, \ldots, N$ {\bf do}\\
\qquad Find a surrogate $g_{k}\in S(f,x_{k-1})$\\
\qquad Find $x_{k}\in \mbox{\rm argmin}_{x\in \O}\,g_{k}(x)$\\
{\bf end for}\\
{\small OUTPUT}: $x_{N}$\\
\hline
\end{tabular}
\end{center}
}

Clearly, the choice of surrogate $g_{k}\in S(f,x_{k-1})$ plays a crucial role in the minimization majorization algorithm.
 In what follows, we consider a particular choice of surrogates for the minimization majorization algorithm; see, e.g., \cite{CZL12,HL,LHY}.  An objective function $f:\O\to \R$ is said to be \emph{majorized} by $\mathcal{M}: \O\times \O\to \R$ if
\begin{equation*}
f(x) \leq \mathcal{M}(x, y)\; \mbox{\rm and } \mathcal{M}(y, y) = f(y) \mbox{ for all }x,y\in\O.
\end{equation*}
Given $x_{k-1}\in\O$, we can define $g_k(x):=\mathcal{M}(x, x_{k-1})$, so that $g_{k}\in S(f,x_{k-1})$. Then the update
\begin{equation*}
x_{k}\in\mbox{\rm arg min}_{x\in \O} \mathcal{M}(x, x_{k-1})
\end{equation*}
defines an minimization majorization algorithm.
As mentioned above, finding an appropriate majorization is an important piece of this algorithm. It has been shown in \cite{CZL12} that the minimization majorization algorithm using distance majorization provides an effective tool for solving many important classes of optimization problems. The key step is to use the following:
\begin{equation*}
d(x; Q)\leq \|x-\Pi(y; Q)\|\;\mbox{\rm and } d(y; Q)=\|y-\Pi(y; Q)\|.
\end{equation*}
In the examples below, we revisit some algorithms based on distance majorization and provide the convergence analysis for these algorithms.

\begin{Example}{\rm Let $\Omega_i $ for $ i = 1,..., m $ be nonempty closed convex subsets of  $\mathbb{R}^n$ such that $\bigcap\limits_{i=1}^m\Omega_i\neq\emptyset$. The problem of finding a point $x^* \in \bigcap\limits_{i=1}^m\Omega_i$ is called the \emph{feasible point problem} for these sets. Consider the problem
\begin{equation}\label{FP1}
\mbox{\rm minimize}\, f(x):=\sum_{i=1}^m[d(x;\Omega_i)]^2, x \in \mathbb{R}^n.
\end{equation}
  With the assumption that $\bigcap\limits_{i=1}^m\Omega_i\neq\emptyset$, $x^*\in \mathbb{R}^n$ is an optimal solution of  \eqref{FP1} if and only if $x^*\in\bigcap\limits_{i=1}^m\Omega_i$. Thus, we only need to consider  \eqref{FP1}.

Let us apply the minimization majorization algorithm for \eqref{FP1}. First, we need to find surrogates for the objective function $f(x)=\sum\limits_{i=1}^m[d(x;\Omega_i)]^2$. Let
$$g_k(x):=\sum_{i=1}^m \|x-\Pi(x_{k-1};\Omega_i)\|^2.$$
Then $g_k\in \mathcal{S}(f,x_{k-1})$ for all $k \in \mathbb{N}$, so the minimization majorization algorithm is given by
$$x_{k}\in \mbox{\rm argmin}_{x\in \R^n}\,g_{k}(x)=\dfrac{1}{m}\sum_{i=1}^m \Pi(x_{k-1}; \O_i).$$
Let
\begin{align*}
h_k(x)&:=g_k(x)-f(x)=\sum_{i=1}^m\left[  \|x-\Pi(x_{k-1};\Omega_i)\|^2-[d(x; \Omega_i)]^2\right].
\end{align*}
We can show that  $h_k(x)$ is differentiable on $\mathbb{R}^n$ and  $\nabla h_k(x)$ is Lipschitz with constant $L=2m$. Moreover, $$h_k(x_{k-1})=\nabla h_k(x_{k-1})=0.$$
The function $g_k(x)-m\|x\|^2$ is convex, so $g_k$ is strongly convex with modulus $\rho=2m$. Using the same notation as in \cite{Mairal}, one has  $g_k\in \mathcal{S}_{L,\rho}(f,x_{k-1})$ with $\rho=L=2m$. By \cite[Proposition~2.8]{Mairal},
$$f(x_k)-V_*\leq \dfrac{m}{k}\|x_0-x^*\|^2\; \text{ for all } k\in \mathbb{N}.$$
}\end{Example}

\begin{Example}{\rm Given a \emph{data set} $S:=\{(a_i, y_i)\}_{i=1}^m$, where $a_i\in \mathbb{R}^p$  and $y_i\in\{-1, 1\}$, consider the \emph{support vector machine} problem
\begin{align*}
&\mbox{\rm minimize}\, \dfrac{1}{2}\|x\|^2\\
&\mbox{\rm subject to } \;y_i\langle a_i,x\rangle\geq 1\; \mbox{\rm for all }i=1,\ldots,m.
\end{align*}
Let $\Omega_i:=\{x\in \mathbb{R}^p \; |\;  y_i\langle a_i,x\rangle\geq 1\} $. Using the quadratic penalty method (see \cite{CZL12}), the support vector machine problem can be solved by the following unconstrained optimization problem:
\begin{equation} \label{svmmm}
\mbox{\rm minimize}\,f(x):=\dfrac{1}{2}\|x\|^2+\dfrac{C}{2}\sum_{i=1}^m[d(x;\Omega_i)]^2, x\in \mathbb{R}^p, C>0.
\end{equation}
Using the minimization majorization algorithm with the surrogates
$$g_k(x)=\dfrac{1}{2}\|x\|^2+\dfrac{C}{2}\sum_{i=1}^m \|x-\Pi(x_{k-1};\Omega_i)\|^2$$ for \eqref{svmmm} yields
$$x_k=\dfrac{C}{1+mC}\sum\limits_{i=1}^m \Pi(x_{k-1},\Omega_i).$$
Let
\begin{align*}
h_k(x):=g_k(x)-f(x)=\dfrac{C}{2}\sum_{i=1}^m\left[  \|x-\Pi(x_{k-1};\Omega_i)\|^2-[d(x;\Omega_i)]^2\right].
\end{align*}
We can show that $\nabla h_k(x)$ is Lipschitz with constant $L=mC$, and $h_k(x_{k-1})=\nabla h_k(x_{k-1})=0.$  Moreover, $g_k$ is strongly convex with parameter $\rho=1+mC$. By \cite[Proposition~2.8]{Mairal}, the minimization majorization method applied for \eqref{svmmm} gives
$$f(x_k)-V_*\leq \dfrac{mC}{2}\left(\dfrac{mC}{mC+2}\right)^{k-1}\|x_0-x^*\|^2 \;\text{ for all }k \in \mathbb{N},$$
where $x_*$ is the optimal solution of \eqref{svmmm} and $V_*$ is the optimal value.
}\end{Example}

In what follows, we apply the minimization majorization algorithm in combination with the log-exponential smoothing technique to solve the smallest intersecting ball problem (\ref{original_prob}). In the first step, we approximate the cost function $\mathcal{D}$ in (\ref{original_prob})  by the log-exponential smoothing function \eqref{smooth_approx}. Then the new function is majorized in order to apply the minimization majorization algorithm. For  $x, y\in \R^n$ and $p>0$, define
$$\mathcal{G}(x, y, p):=p\ln\sum_{i=1}^m\exp\left(\frac{\sqrt{\|x-\Pi(y; \O_i)\|^2+p^2}}{p}\right).$$
Then $\mathcal{G}(x, y, p)$ serves as a majorization of the log-exponential smoothing function (\ref{smooth_approx}). From Proposition~\ref{lemma2}, for $p>0$ and $y\in \R^n$, the function $\mathcal{G}(x, y, p)$ with variable $x$ is strongly convex on any bounded set and continuously differentiable with Lipschitz gradient on $\R^n$.

Our algorithm is explained as follows. Choose a small number $\bar{p}$. In order to solve the smallest intersecting ball problem (\ref{original_prob}), we minimize its log-exponential smoothing approximation \eqref{smooth_approx}:
\begin{equation}\label{EXP}
\mbox{\rm minimize}\; \mathcal{D}(x, \bar p)\; \; \mbox{\rm subject to } \; x\in \O.
\end{equation}
Pick an initial point $x_0\in \O$ and apply the minimization majorization algorithm with
\begin{align}\label{MM}
x_{k}:=\mbox{\rm arg}\min_{x\in \O}\mathcal{G}(x, x_{k-1}, \bar p).
\end{align}
The algorithm is summarized by the following.\\[1ex]
{\small{\bf Algorithm 2}.
\begin{center}
\begin{tabular}{| l |}
\hline
{\small INPUT}: $\O$, $\bar p>0$, $x_0 \in \O$, $m$ target sets $\O_i$, $i=1, \ldots, m$, $N$\\
{\bf for } $k=1, \ldots, N$ {\bf do}\\
\qquad {\bf use a fast gradient algorithm to solve approximately} \\
\qquad \qquad $x_k:= \mbox{\rm arg}\min_{x\in \O}\mathcal{G}(x, x_{k-1}, \bar p)$\\
{\bf end for}\\
{\small OUTPUT}: $x_{N}$\\
\hline
\end{tabular}
\end{center}
}

\begin{Proposition}\label{limit_point}
Given $\bar p>0$ and $x_0\in \O$, the sequence $\{x_k\}$ of exact solutions $x_k:= \mbox{\rm arg}\min_{x\in \O}\mathcal{G}(x, x_{k-1}, \bar p)$ generated by Algorithm 2 has a convergent subsequence.
\end{Proposition}
{\it Proof} Denoting $\alpha:= \mathcal{D}(x_0, \bar p)$ and using Theorem \ref{lm3.1}(v) imply that the level set
$$\mathcal{L}_{\leq \alpha}:=\{x\in \O \ | \ \mathcal{D}(x, \bar p) \leq \alpha\}$$
is bounded. For any $k\geq 1$, because $\mathcal{G}(\cdot, x_{k-1},\bar p)$ is a surrogate of $\mathcal{D}(\cdot,\bar p)$ at $x_{k-1}$, one has
\begin{align*}
\mathcal{D}(x_k,\bar p) \leq \mathcal{G}(x_k, x_{k-1},\bar p) & \leq \mathcal{G}(x_{k-1},  x_{k-1},\bar p) = \mathcal{D}(x_{k-1},\bar p).
\end{align*}
It follows that
\begin{align*}
\mathcal{D}(x_k,\bar p) \leq \mathcal{D}(x_{k-1},\bar p)\leq \ldots \leq  \mathcal{D}(x_{1},\bar p)\leq  \mathcal{D}(x_{0}, \bar p).
\end{align*}
This implies $\{x_k\} \subset  \mathcal{L}_{\leq \alpha}$ which is a bounded set, so $\{x_k\}$ has a convergent subsequence. $\h$

The convergence of  the minimization majorization algorithm depends on the algorithm map
\begin{equation}\label{AM}
\psi(x):=\mbox{\rm arg}\min_{y \in \O}\mathcal{G}(y, x, \bar p) = \mbox{\rm arg}\min_{y\in \O}\bigg\{\bar p\ln\sum_{i=1}^m\exp\left(\frac{\sqrt{
\|y-\Pi(x; \O_i)\|^2+{\bar p}^2}}{\bar p}\right)\bigg\}.
\end{equation}
In the theorem below, we show that the conditions in \cite[Proposition~1]{CZL12} are satisfied.
\begin{Theorem}\label{CV1} Given $\bar p>0$, the function $\mathcal{D}(\cdot, \bar p)$ and the  algorithm map $\psi: \O\to \O$ defined by {\rm (\ref{AM})} satisfy the following conditions:\\[1ex]
{\rm\bf (i)} For any $x_0\in \O$, the level set
\begin{equation*}
\mathcal{L}(x_0):=\{x\in \O\; |\; \mathcal{D}(x,\bar p)\leq \mathcal{D}(x_0, \bar p)\}
\end{equation*}
is compact.\\
{\rm\bf (ii)} $\psi$ is continuous on $\O$.\\
{\rm\bf (iii)} $\mathcal{D}(\psi(x), \bar p)<\mathcal{D}(x, \bar p)$ whenever $x\neq \psi(x)$.\\
{\rm\bf (iv)} Any fixed point $\ox$ of $\psi$ is a minimizer of $\mathcal{D}(\cdot, \bar p)$ on $\O$.
\end{Theorem}
{\it Proof} Observe that the function $\mathcal{D}(\cdot,\bar p)$ is continuous on $\O$. Then the level set $\mathcal{L}(x_0)$ is compact for any initial point $x_0$ since $\mathcal{D}(\cdot, \bar p)$ is coercive by Theorem~\ref{lm3.1}(v), and hence {\bf(i)} is satisfied. From the strict convexity on $\O$ of $\mathcal{G}(\cdot, x, \bar p)$ guaranteed by Proposition~\ref{lemma2}, we can show that the algorithm map $\psi: \O\to \O$ is a single-valued mapping. Let us prove that $\psi$ is continuous. Take an arbitrary sequence $\{x_k\} \subset \O$, $x_k \to \ox\in \O$ as $k\to \infty$. It suffices to show that the sequence $y_k:=\psi(x_k)$ tends to $\psi(\ox)$. It follows from the continuity of $\mathcal{D}(\cdot, \bar p)$ that $\mathcal{D}(x_k, \bar p) \to \mathcal{D}(\ox, \bar p)$, and hence we can assume
$$\mathcal{D}(x_k, \bar p) \leq \mathcal{D}(\ox, \bar p) +\delta,$$
for all $k\in \N$, where $\delta$ is a positive constant. One has the estimates
\begin{equation*}
\mathcal{D}(\psi(x_k), \bar p) \leq \mathcal{G}(\psi(x_k), x_k, \bar p) \leq \mathcal{G}(x_k, x_k, \bar p) =\mathcal{D}(x_k, \bar p) \leq \mathcal{D}(\ox, \bar p) +\delta,
\end{equation*}
which imply that $\{y_k\}$ is bounded by the coerciveness of $\mathcal{D}(\cdot, \bar p)$. Consider any convergent subsequence $\{y_{k_\ell}\}$ with the limit $z$. Since $y_{k_\ell}$ is a solution of the smooth optimization problem $\min_{y\in \O}\mathcal{G}(y, x_{k_\ell}, \bar p)$, by the necessary and sufficient optimality condition \eqref{optimality} for the given smooth convex constrained optimization problem, we have
\begin{equation*}
\la \nabla \mathcal{G}(y_{k_\ell}, x_{k_\ell}, \bar p), x-y_{k_\ell} \ra \geq 0 \; \; \mbox{\rm for all}\;x\in \O.
\end{equation*}
This is equivalent to
\begin{equation*}
\bigg \la \sum_{i=1}^m \frac{\lambda_i(y_{k_\ell}, \bar p)}{g_i(y_{k_\ell}, \bar p)}\left(y_{k_\ell}- \Pi(x_{k_\ell}; \O_i)\right), x-y_{k_\ell} \bigg \ra \geq 0 \; \; \mbox{\rm for all}\;x\in \O,
\end{equation*}
where
\begin{equation*}
g_i(y_{k_\ell}, \bar p)=\sqrt{\|y_{k_\ell}-\Pi(x_{k_\ell}, \O_i)\|^2 + {\bar p}^2} \; \mbox{\rm and }
\lambda_i(y_{k_\ell}, \bar p)=\frac{\exp\left(g_i(y_{k_\ell}, \bar p)/\bar p\right)}{\sum_{i=1}^m \exp\left(g_i(y_{k_\ell}, \bar p)/\bar p\right)}.
\end{equation*}
Since the Euclidean projection mapping to a nonempty closed convex set is continuous, by passing to a limit, we have
\begin{equation*}
\la \nabla \mathcal{G}(z, \ox, \bar p), x-z\ra \geq 0 \; \; \mbox{\rm for all}\;x\in \O.
\end{equation*}
Thus, applying \eqref{optimality} again implies that $z$ is also an optimal solution of the problem $\min_{y\in \O}\mathcal{G}(y,\ox, \bar p)$. By the uniqueness of solution and $\psi(\ox)=\mbox{\rm arg}\min_{y\in \O}\mathcal{G}(y,\ox, \bar p)$, one has that $z=\psi(\ox)$ and $\{y_{k_\ell}\}$ converges to $\psi(\ox)$. Since this conclusion holds for all convergent subsequences of the bounded sequence $\{y_k\}$, the sequence $\{y_k\}$ itself converges to $\psi(\ox)$, which shows that {\bf(ii)} is  satisfied. Let us verify that $\mathcal{D}(\psi(x), \bar p)< \mathcal{D}(x, \bar p)$ whenever $\psi(x)\neq x$. Observe that $\psi(x)=x$ if and only if $\mathcal{G}(x, x, \bar p)=\min_{y\in \O} \mathcal{G}(y, x, \bar p)$. Since $\mathcal{G}(y, x, \bar p)$ has a unique minimizer, we have the strict inequality $\mathcal{G}(\psi(x), x, \bar p)<\mathcal{G}(x, x, \bar p)$ whenever $x$ is not a fix point of $\psi$. Combining with $\mathcal{D}(\psi(x), \bar p) \leq \mathcal{G}(\psi(x), x, \bar p)$ and $\mathcal{D}(x, \bar p) = \mathcal{G}(x, x, \bar p)$, we arrive at the conclusion {\bf(iii)}.

Finally, we show that, any fixed point $\ox$ of algorithm map $\psi(x)$ is a minimizer of $\mathcal{D}(x,\bar p)$ on $\O$. Fix any $\ox\in \O$ such that $\psi(\ox)=\ox$. Then $\mathcal{G}(\ox, \ox, \bar p)=\min_{y\in \O}\mathcal{G}(y, \ox, \bar p)$, which is equivalent to
\begin{equation*}
\la \nabla \mathcal{G}(\ox, \ox, \bar p), x-\ox \ra \geq 0 \; \mbox{\rm for all}\;x\in \O.
\end{equation*}
This means
\begin{equation*}
\bigg \la \sum_{i=1}^m \frac{\lambda_i(\ox, \bar p)}{g_i(\ox, \bar p)}\left(\ox- \Pi(\ox; \O_i)\right), x-\ox \bigg \ra \geq 0 \; \mbox{\rm for all}\;x\in \O,
\end{equation*}
where
\begin{equation*}
g_i(\ox, \bar p)=\sqrt{\|\ox-\Pi(\ox, \O_i)\|^2 + {\bar p}^2}=\sqrt{d(\ox; \O_i)^2+{\bar p}^2}=G_i(\ox, \bar p)
\end{equation*}
and
\begin{equation*}
\lambda_i(\ox, \bar p)=\frac{\exp\left(g_i(\ox, \bar p)/\bar p\right)}{\sum_{i=1}^m \exp\left(g_i(\ox, \bar p)/\bar p\right)}=\Lambda_i(\ox, \bar p).
\end{equation*}
This inequality, however, is equivalent to the inequality $\la \nabla \mathcal{D}(\ox, \bar p), x-\ox \ra \geq 0$ for all $x\in \O$, which in turn holds if and only if $\ox$ is a minimizer of $\mathcal{D}(x,\bar p)$ on $\O$.
$\h$

\begin{Corollary}\label{convergence cor} Given $\bar p>0$ and $x_0\in \O$, the sequence $\{x_k\}$ of exact solution $x_k:= \mbox{\rm arg}\min_{x\in \O}\mathcal{G}(x, x_{k-1}, \bar p)$ generated by Algorithm 2 has a subsequence that converges to an optimal solution of {\rm (\ref{EXP})}. If we suppose further that problem {\rm (\ref{EXP})} has a unique optimal solution, then $\{x_k\}$ converges to this optimal solution.
\end{Corollary}
{\it Proof} It follows from Proposition \ref{limit_point} that $\{x_k\}$ has a subsequence $\{x_{k_\ell}\}$ that converges to $\ox$. Applying \cite[Proposition~1]{CZL12} implies that $\|x_{k_\ell+1}-x_{k_\ell}\|\to 0$ as $k\to \infty$. From the continuity of the algorithm map $\psi$ and the equation $x_{k_\ell+1}=\psi(x_{k_\ell})$, one has that $\psi(\ox)=\ox$. By Theorem \ref{CV1}(iv), the element $\ox$ is an optimal solution of {\rm (\ref{EXP})}. The last conclusion is obvious.
$\h$

In what follows, we apply \emph{Nesterov's accelerated gradient method} introduced in \cite{n,n83} to solve (\ref{MM}) approximately. Let $f: \R^n\to \R$ be a a smooth convex function
with Lipschitz gradient. That is, there exists $\ell\geq 0$ such that
\begin{equation*}
\|\nabla f(x)-\nabla f(y)\|\leq \ell \|x-y\|\; \mbox{ for all }x, y\in \R^n.
\end{equation*}
Let $\O$ be a nonempty closed convex set. In his seminal papers \cite{n,n83}, Nesterov considered the optimization problem
\begin{equation*}
\mbox{\rm minimize } f(x)\; \mbox{\rm subject to }x\in \O.
\end{equation*}
For $x\in \R^n$, define
\begin{equation*}
T_\O(x):=\mbox{\rm arg min }\{\la \nabla f(x), y-x\ra +\dfrac{\ell}{2}\|x-y\|^2\; |\; y\in \O\}.
\end{equation*}
Let $d: \R^n\to \R$ be a strongly convex function with parameter $\sigma>0$. Let $x_0\in \R^n$  such that
$$x_0=\mbox{\rm arg min }\{d(x)\; |\; x\in \O\}.$$
Further, assume that $d(x_0)=0$.

For simplicity, we choose $d(x)=\frac{1}{2}\|x-x_0\|^2$, where $x_0\in \O$, so $\sigma=1$. It is not hard to see that
\begin{equation*}
y_k=T_\O(x_k)=\Pi(x_k-\dfrac{\nabla f(x_k)}{\ell}; \O).
\end{equation*}
Moreover,
\begin{equation*}
z_k=\Pi(x_0-\dfrac{1}{\ell}\sum_{i=0}^k\dfrac{i+1}{2}\nabla f(x_i);\O).
\end{equation*}

Nesterov's accelerated gradient algorithm is outlined as follows.\\[1ex]
{\small{\bf Algorithm 3}.
\begin{center}
\begin{tabular}{| l |}
\hline
{\small INPUT}: $f$, $\ell$, $x_0\in \O$\\
set $k = 0$\\
{\bf repeat }\\
\qquad find $y_k := T_\O(x_k)$\\
\qquad find $z_k :=\mbox{\rm arg min}\big\{\dfrac{\ell}{\sigma}d(x)+\sum_{i=0}^k\dfrac{i+1}{2}[f(x_i)+\la \nabla f(x_i), x-x_i\ra]\; \big|\; x\in \O\big\}$\\
\qquad set $x_k :=\dfrac{2}{k+3}z_k+\dfrac{k+1}{k+3}y_k$\\
\qquad set $k := k+1$\\
{\bf until a stopping criterion is satisfied.}\\
{\small OUTPUT}: $y_k$.\\
\hline
\end{tabular}
\end{center}
}

It has been experimentally observed that the algorithm is more effective if, instead of choosing a small value $p$ ahead of time, we change its value using an initial value $p_0$ and define $p_{s}:=\sigma p_{s-1}$, where $\sigma\in (0,1)$. \\[1ex]
{\small{\bf Algorithm 4}.
\begin{center}
\begin{tabular}{| l |}
\hline
{\small INPUT}: $\O$, $\epsilon>0$, $p_0>0$, $x_0 \in \O$, $m$ target sets $\O_i$, $i=1, \ldots, m$, $N$\\
set $p=p_0$, $y=x_0$\\
{\bf for } $k=1, \ldots, N$ {\bf do}\\
\qquad {\bf use Nesterov's accelerated gradient method to solve approximately}\\
\qquad \qquad $y := \mbox{\rm arg}\min_{x\in \O}\mathcal{G}(x, y, p)$\\
\qquad set $p := \sigma p$\\
{\bf end for}\\
{\small OUTPUT}: $y$\\
\hline
\end{tabular}
\end{center}
}
\begin{Remark}\label{rmmm} {\rm {\bf (i)} When implementing this algorithm, we usually desire to maintain $p_k>\epsilon$, where $\epsilon<p_0$. So the factor $\sigma$ can be calculated based on the desired number of iterations $N$ to be run, i.e., $\sigma=(\epsilon/p_0)^{1/N}$.\\[1ex]
{\bf (ii)} In Nesterov's accelerated gradient method, at iteration $k$, we often use the stopping criterion
$$\|\nabla_x \mathcal{G}(x, y, p)\|<\gamma_k,$$
where $\gamma_0$ is chosen and $\gamma_k=\tilde{\sigma}\gamma_{k-1}$ for some $\tilde{\sigma}\in (0,1)$. The factor $\tilde{\sigma}$ can be calculated based on the number of iterations $N$ and the lower bound $\tilde{\epsilon}>0$ for $\gamma_k$ as above.}

\end{Remark}

\section{Numerical Implementation}
We implement Algorithm 4 to solve the generalized Sylvester problem in a number of  examples. In each of the following examples, we implement Algorithm 4 with the following parameters described in this algorithm and Remark \ref{rmmm}: $\epsilon = 10^{-6}, \tilde{\epsilon}=10^{-5}, p_0=5, \gamma_0= .5,$ and $N=10$. Observations suggest that when the number of dimensions is large, speed is improved by starting with the relatively high $\gamma_0$ and $p_0$ and decreasing each (thereby reducing error) with each iteration. Choosing $\sigma, \tilde{\sigma}$ as described in Remark \ref{rmmm} ensures that the final iterations are of desired accuracy.
The approximate radii of the smallest intersecting ball that corresponds to the approximate optimal solution $x_k$ is $r_k:=\mathcal{D}(x_k)$.
%\begin{figure}[hbt]
%\hspace{2cm} \includegraphics[width=12cm]{Sylvester_for_disks.eps}
%\caption{A generalized Sylvester problem for disks in $\R^2$}
%\label{fig:1}
%\end{figure}
\begin{figure}[!ht]
\begin{center}
 \includegraphics[width=5cm]{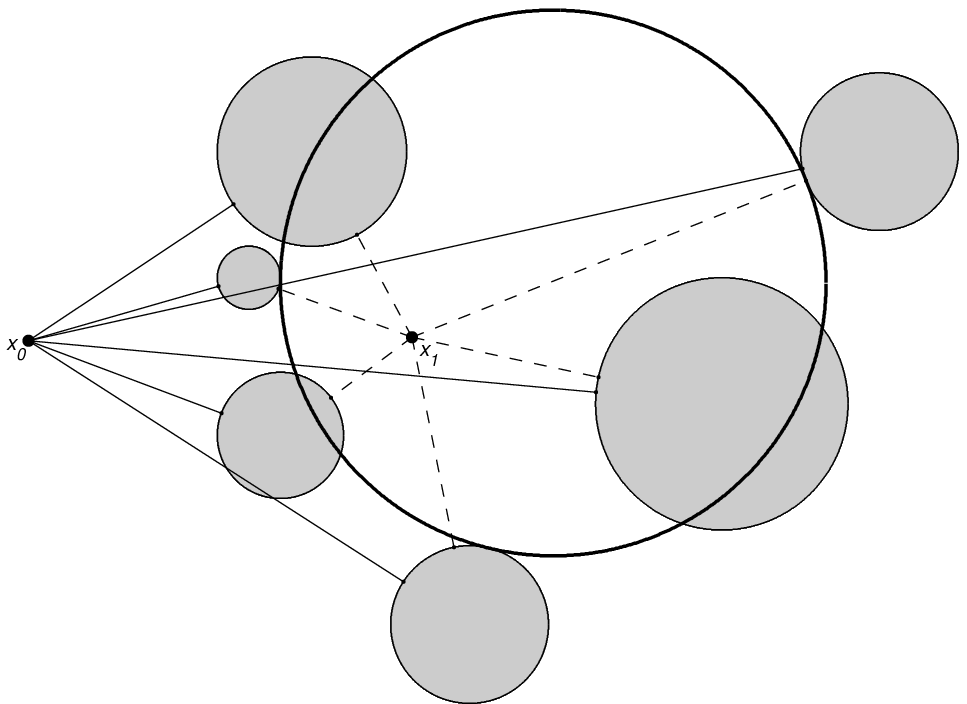}
\caption{An illustration of minimization majorization algorithm. }
\label{fig:1}
\end{center}
\end{figure}
\begin{Example}
{\rm Let us first apply Algorithm 4 to an unconstraint generalized Sylvester problem (\ref{original_prob}) in $\R^2$ in which $\O_i$ for $i=1, \ldots, 6$ are disks with centers at $(-6, 9)$, $(12, 9)$, $(-1,-6)$, $(-8, 5)$, $(-7, 0)$, $(7,1)$ with radii $3, 2.5, 2.5, 1,  2,  4$, respectively. This setup is depicted in Figure \ref{fig:1}. A simple MATLAB program yields an approximate smallest intersecting ball with center $x_* \approx  (1.65,4.83)$ with the approximate radius $r_*\approx 8.65$. Figure \ref{fig:1} shows a significant move toward the optimal solution of the problem in one step of the minimization majorization algorithm.
}
\end{Example}

\begin{figure}[!ht]
\begin{center}
 \includegraphics[width=5cm]{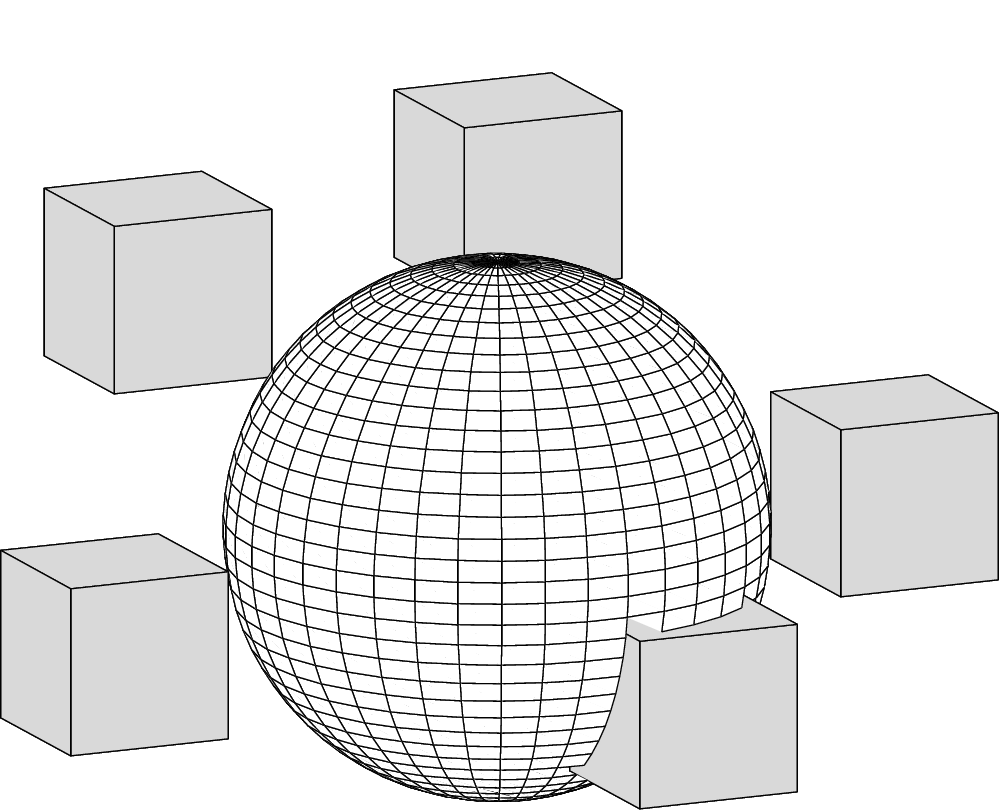}

\caption{A smallest intersecting ball problem for cubes in $\R^3$.}
\label{fig:2}
\end{center}
\end{figure}

\begin{Example}\label{ex2}{\rm We consider the smallest intersecting ball problem in which the target sets are square boxes in $\R^n$.
In $\R^n$, a square box $S(\omega, r)$ with center $\omega=(\omega_1, \ldots, \omega_n)$ and radius $r$ is the set
$$S(\omega, r):=\{x=(x_1, \ldots, x_n)\ |\  \max\{|x_1-\omega_1|, \ldots, |x_n-\omega_n|\}\leq r\}.$$
Note that the Euclidean projection from $x$ to $S(\omega, r)$ can be expressed componentwise as follows
$$[\Pi(x; S)]_i=\begin{cases}
\omega_i-r_i&\mbox{if } x_i +r \leq \omega_i,\\
x_i&\mbox{if } \omega_i -r \leq x_i \leq \omega_i +r,\\
\omega_i+r&\mbox{if } \omega_i +r \leq x_i.
\end{cases}$$
Consider the case where $n=3$. The target sets are $5$ square boxes with centers $(-5, 0, 0)$, $(1, 4, 4)$,
$(0, 5, 0)$, $(-4, -3, 2)$ and $(0, 0, 5)$ and radii $r_i = 1$ for $i =1,\ldots, 5$. Our results show that both Algorithm 4 and the subgradient method give an approximate SIB radius $r_*\approx 3.18$; see Figure \ref{fig:2}.
}
\end{Example}

\begin{center}
\vspace{15pt}
\begin{figure}[h]
\begin{minipage}{2.3in}
   \includegraphics[width=3.0in]{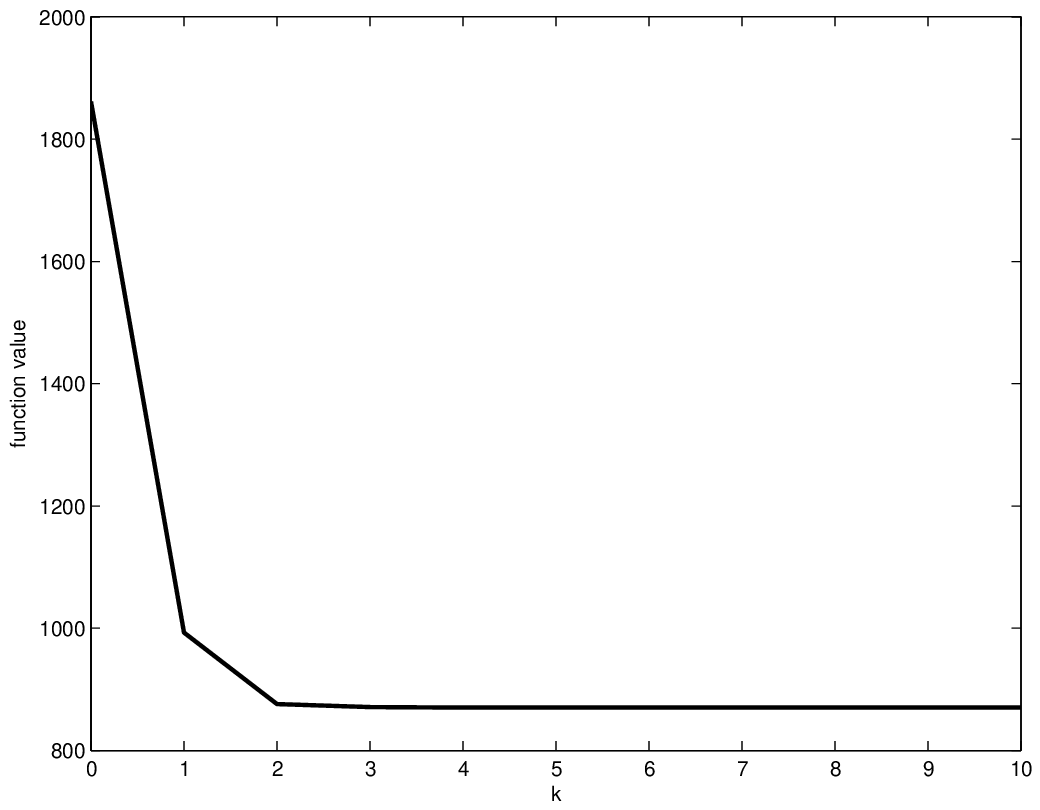}\\
\end{minipage}
~\hfill~
\begin{minipage}[t]{0.6\textwidth}
\hspace{90pt}\begin{tabular}{|c|c|}
\hline
\multicolumn{2}{|c|}{MATLAB RESULT} \\
\hline
$k$ &  $r_k$ \\
\hline
 0 &  1861.36441\\
   1 &  992.32230\\
   2 &  875.27618\\
   3 &  870.47714\\
   4 &  869.94621\\
   5 &  869.82005\\
   6 &  869.79982\\
   7 &  869.79676\\
   8 &  869.79628\\
   9 &  869.79621\\
   10&  869.79619\\
\hline
\end{tabular}
\end{minipage}
\vspace{-15pt}
\caption{\small A smallest intersecting ball problem for cubes in high dimensions.}
\end{figure}
\label{fig:3}
\end{center}

\begin{Example}{\rm
Now we illustrate the performance of Algorithm 4 in high dimensions with the same setting in Example \ref{ex2}. We choose a modification of the pseudo random sequence from \cite{ZTS05} with $a_0=7$ and for $i=1, 2, \ldots$
$$a_{i+1}=  \textrm {mod}(445a_i + 1, 4096) ,\;\; b_i=\frac{a_i}{40.96}.$$
The radii $r_i$ and centers $c_i$ of the square boxes are successively set to $b_1, b_2, \ldots$ in the following order
$$10r_1, c_1(1), \ldots, c_1(n); 10r_2, c_2(1), \ldots, c_2(n); \ldots; 10r_m, c_m(1), \ldots, c_m(n).$$
Consider $m=100, n=1000$. Figure 3 shows the approximate values of the radii $r_k$ for $k=0, \ldots, 10$.

\begin{figure}[!ht]
\begin{center}
 \includegraphics[width=7cm]{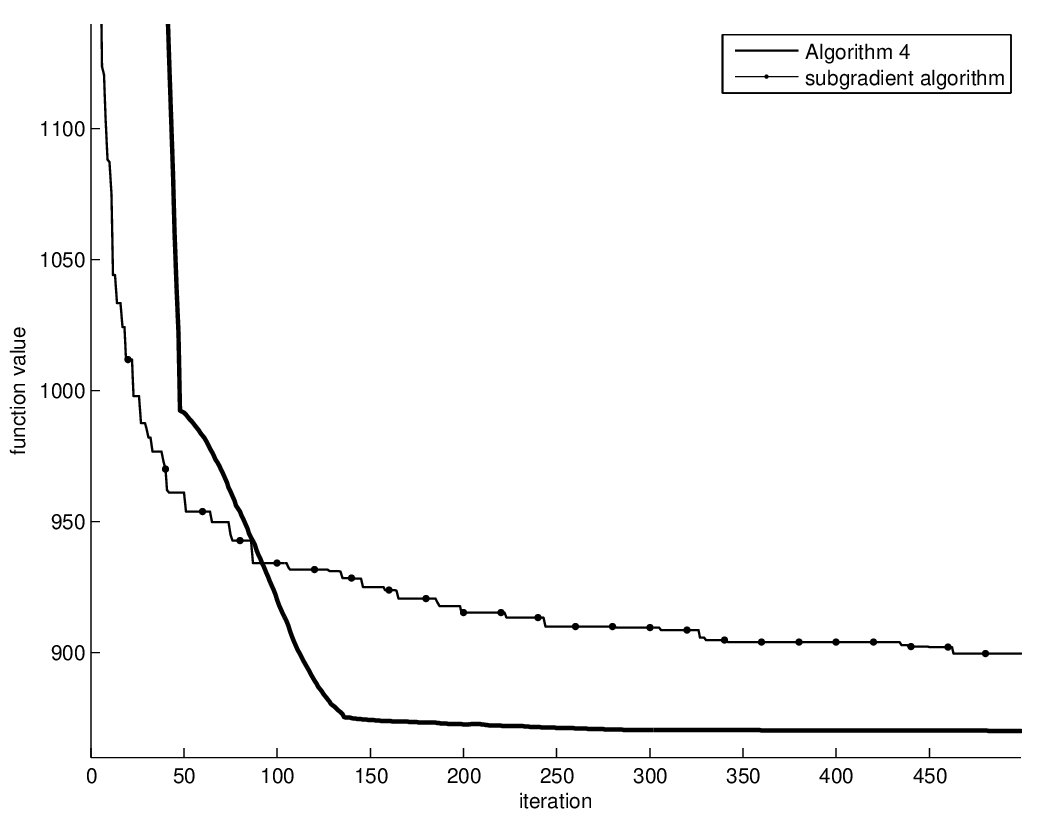}
\caption{Comparison between minimization majorization algorithm and subgradient algorithm.}
\label{fig:4}
\end{center}
\end{figure}

We also implement algorithm 4 in comparison with the subgradient algorithm. From our numerical results, we see that Algorithm 4 performs much better than the subgradient algorithm in both accuracy and speed. In the case where the number of  target sets is large or the dimension $m$ is high, the subgradient algorithm seems to be stagnated but Algorithm 4 still performs well. Figure 4 shows that comparison between Algorithm 4 and the subgradient algorithm. Note that in Algorithm 4, we count every iteration of Nesterov's accelerated gradient method in the total iteration count along the horizontal axis.  Thus the ``sharp corner" that can be seen at 50 iterations represents the transition form $x_0$ to $x_1$, and a subsequent recalculation of $p$ by $p=\sigma p$.
}
\end{Example}

\section{Concluding Remarks}
\label{s:Concluding}
Based on the log-exponential smoothing technique and the minimization majorization algorithm, we have developed an effective numerical algorithm for solving the smallest intersecting ball problem. The problem under consideration %which involves minimizing functions of distance functions to convex sets%
 not only generalizes the Sylvester smallest enclosing circle problem, but also opens up the possibility of applications to other problems of constrained optimization, especially those that appear frequently in machine learning. Our numerical examples show that the algorithm works well for the problem in high dimensions. %Although we have proved a number of convergence results, more careful study on the convergence rate of the algorithm will be considered in our future research.%
 Although a number of key convergence results are contained in this paper, our future work will further develop an understanding of the convergence rate of this algorithm

{\bf Acknowledgement.} The author would like to thank Prof. Jie Sun for giving comments that help improve the presentation of the paper.

\end{document}